\documentclass[letterpaper,10pt,conference]{ieeeconf}
\IEEEoverridecommandlockouts
\IEEEoverridecommandlockouts                              

\makeatletter

\let\NAT@parse\undefined

\makeatother

\usepackage{cite}
\usepackage{hyperref,url,cite}
\usepackage{amsmath,amssymb,amsfonts}
\usepackage{graphicx}
\usepackage{textcomp}
\usepackage{xcolor}
\def\BibTeX{{\rm B\kern-.05em{\sc i\kern-.025em b}\kern-.08em
    T\kern-.1667em\lower.7ex\hbox{E}\kern-.125emX}}
\newtheorem{theorem}{Theorem}

\newtheorem{proposition}{Proposition}

\title{\LARGE \bf 
Real-time control of metro train dynamics with minimization of train time-headway variance}

\author{Florian Schanzenb\"acher$^{1*}$, Nadir Farhi$^{2}$, Fabien Leurent$^{3}$, G\'erard Gabriel$^{4}$
\thanks{$^{1}$ RATP, Paris, France and Universit\'e Paris-Est, France.}%
\thanks{$^{2}$ Universit\'e Paris-Est, Ifsttar/Cosys/Grettia, F-77447 France.}
\thanks{$^{3}$ LVMT, ENPC, France. } %
\thanks{$^{4}$ RATP, Paris, France.}%
\thanks{* Corresponding author. {\texttt{florian.schanzenbacher@ratp.fr}}} }

\begin{document}
\maketitle
\thispagestyle{empty}
\pagestyle{empty}
\begin{abstract}
We present here a real-time control model for the train dynamics in a linear metro line system.
The model describes the train dynamics taking into account average passenger arrival rates on platforms,
including control laws for train dwell and run times, based on the feedback of the train dynamics.
The model extends a recently developed Max-plus linear traffic model with demand-dependent dwell times and a run time control. 
The extension permits the elimination of eventual irregularities on the train time-headway.
The resulting train dynamics are interpreted as a dynamic programming system of a stochastic optimal control problem of a Markov chain.
The train dynamics still admit a stable stationary regime with a unique average growth rate interpreted as the asymptotic
average train time-headway. Moreover, beyond the transient regime of the train dynamics, our extension guarantees uniformity in time of the train time-headways at every platform.
\end{abstract}


\section{Introduction and literature review}
This paper is part of a series of works being realized as a research project at RATP, the Paris metro operator. 
The research project aims to develop a passenger demand-dependent traffic management and real-time control algorithm for high density metro lines. The characteristic of the algorithm is its sensibility to passenger demand, i.e. a premier aim is to control the traffic while maintaining a high frequency, satisfying the passenger demand.

The modeling approach is based on former works from~\cite{Bres91}, where a traffic control for metro loop lines has been proposed. This approach has served as a basis for the authors of~\cite{Fer05} which have presented a model predictive control algorithm for metro loop lines aiming for headway regularity one the one side, and schedule adherence, on the other side. Lately, the authors of~\cite{SCF16} have applied this model predictive control algorithm to a part of a mass transit railway line in the area of Paris.
Still on the same basis, the authors of~\cite{LDYG17a,LYG17b} have proposed a traffic control algorithm for high density metro lines, taking into account the passenger demand.

Starting with a recently developed traffic model for metro lines, with static dwell and run times, in~\cite{FNHL16}, the authors of~\cite{SFCLG17} have extended this approach to a metro line with a junction. Applied to a RATP metro line in Paris, this traffic model has precisely reproduced the microscopic characteristics of the line, e.g. the bottleneck segment of the line, the average frequency depending on dwell, run, safe separation times and the number of segments and trains.
In~\cite{SFLG18a}, the authors of this paper here, have replaced the static dwell and run times by a dynamic passenger demand-dependent dwell time. A run time control has been introduced, avoiding the amplification of perturbations, i.e. a cascade effect.
These traffic models have been shown to be linear in Max-plus algebra, which allows to derive analytic formulas for the average headway and frequency on the line.

In this paper, we refer to a Max-plus Theorem, see~\cite{BCOQ92, CCGMQ98, Gov07}, having already been used in our preceding works. 
In the following, we present an enhanced version of the traffic model which we have proposed in~\cite{SFLG18a}. It keeps its advantages, i.e. the passenger demand-dependency and the cancellation of a possible cascade effect.
As a new element, the here presented model includes a train position control. This means, we modify the model of~\cite{SFLG18a} such that a minimization of train time-headway variances over the metro line is realized, while maintaining the average frequency on the line, which is a crucial issue for mega cities as Paris, with regard to the very high passenger demand.
In the upcoming sections we recall first the modeling approach this paper is based on, before presenting the new traffic model with headway variance minimization. We underpin the practical relevance of this new traffic model with simulation results on a RATP metro line and conclude with an outlook on upcoming work.

\section{Review of the modeling approach}

\begin{figure}[htbp]
  \centering
  \includegraphics[width=0.5\textwidth]{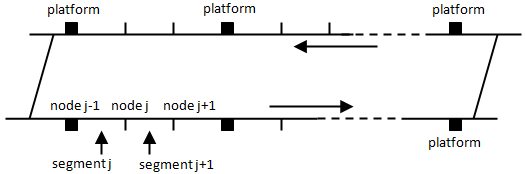} 
  \caption {Schema of a metro loop line and the corresponding notation.}
  \label{fig_line}
\end{figure}

We consider a linear metro line without junction, as in Fig.~\ref{fig_line}.
The line is discretized in space into a number of segments where the length of every segment
is bigger than the train length, as in~\cite{FNHL16}.
The train dynamics is modeled here by applying run and dwell time control laws under train speed and safe separation constraints.
The passenger travel demand is modeled with a static model, considering average arrival rates to every platform.

\subsection{Notations}


We use the following notations, similar to the ones used in~\cite{SFLG18a}.
For the train dynamics:
\\~\\
\begin{tabular}{ll}
    $n$ & the number of segments of the linear metro line.\\
    $m$ & the number of running trains on the metro line.
\end{tabular}
\begin{tabular}{ll}
  $d_{j}^{k}$ & the $k^{th}$ departure time of trains on segment $j$. \\
  $a_{j}^{k}$ & the $k^{th}$ arrival time of trains on segment $j$.\\
  $r^k_{j}$ & the $k^{th}$ running time on segment $j$ \\
      & (between nodes $j-1$ and $j$). \\
  $w^k_{j}$ & $=d^k_{j}-a^k_{j}$ the $k^{\text{th}}$ dwell time on node $j$.
\end{tabular}  
\begin{tabular}{ll}
  $t^k_{j}$ & $=r^k_{j}+w^k_{j}$ the $k^{\text{th}}$ travel time from node \\
		& $j-1$ to node $j$.\\
  $g^k_{j}$ & $=a^k_{j}-d^{k-1}_{j}$ the $k^{\text{th}}$ safe separation time \\
                & (or close-in time) at node $j$.		
\end{tabular}  
\begin{tabular}{ll}		
  $h^k_{j}$ & $=d^k_{j}-d^{k-1}_{j} = g^k_{j}+w^{k}_{j}$ the $k^{\text{th}}$\\
      & departure time-headway at node $j$.\\
  $s^k_{j}$ & $=g^{k+b_{j}}_{j}-r_{j}$.\\    
  $b_{j}$ & $\in \{0,1\}$. It is $0$ (resp. $1$) if there is no train\\
              & (resp. one train) on segment $j$ at time zero.\\
  $\bar{b}_{j}$ & $= 1 - b_{j}$.
\end{tabular}~~\\~~

The average on $k$ and $j$ of the quantities above are denoted 
$r, w, t, g, h$ and $s$. 
We have the following relationships~\cite{FNHL16}.
\begin{align}
    & g = r + s, \label{form1} \\
    & t = r + w, \label{form2} \\
    & h = g + w = r + w + s = t + s. \label{form3}
\end{align} 

For the passenger demand:
\\~\\
\begin{tabular}{ll}  
  $\lambda_{ij}$ & passenger travel demand from platform~$i$ to \\
      & platform~$j$, when $i$ and $j$ denote platforms, and\\
      & $\lambda_{ij} = 0$ if $i$ or $j$ is not a platform node.\\
  $\lambda^{\text{in}}_i$ & $=\sum_j \lambda_{ij}$ the average passenger arrival rate on \\
      & origin platform $i$ to any destination platform. 
\end{tabular}  
\begin{tabular}{ll}      
  $\lambda^{\text{out}}_j$ & $=\sum_i \lambda_{ij}$ the average passenger arrival rate on \\
      & any origin platform to destination platform $j$. \\
  $\alpha^{\text{in}}_j$ & average passenger upload rate on platform $j$.  \\  
  $\alpha^{\text{out}}_j$ & average passenger download rate on platform $j$.
\end{tabular}  
\begin{tabular}{ll}    
  $\frac{\sum_{i} \lambda_{ij} h_{i}}{ \alpha_{j}^{out}}$ & time for passenger download at\\
  & platform $j$. \\
  $\frac{\lambda^{\text{in}}_j}{ \alpha^{\text{in}}_j} \; h_j$ & time for passenger upload at platform $j$. \\
\end{tabular}  
\\~\\

We approximate here $\sum_{i} \lambda_{ij} h_{i} / \alpha_{j}^{out}$ by $(\lambda^{\text{out}}_j / \alpha^{\text{out}}_j)\; h_j$.
We then define $x_j$ as a passenger demand parameter
\begin{equation}
   x_j = \left( \frac{\lambda^{\text{out}}_j}{\alpha^{\text{out}}_j} + \frac{\lambda^{\text{in}}_j}{\alpha^{\text{in}}_j} \right),
\label{eq-x}   
\end{equation}
such that $x_j h_j$ represents the time needed for passenger alighting and boarding at platform $j$.
We notice that capacity limits for trains and platforms are not modeled here.
Furthermore, we suppose that the passenger upload and download rate are significantly larger than the travel demand, such that
$x_{j} << 1$. We also define
\begin{equation}
  X_{j} = \frac{x_{j}}{1 - x_{j}},
\label{eq-X}
\end{equation}
in such a way that $x_{j} h_{j} = X_{j} g_{j}$. 

Lower bounds are written
$\underline{r}_j, \underline{w}_j, \underline{t}_j, \underline{g}_j, \underline{s}_j, \underline{h}_j$.
The minimum run and dwell time $\underline{r}_j, \underline{w}_j$ are directly given by the infrastructure constraints of the metro system, as well as $\underline{g}_j$, i.e. the minimum time span between the $(k-1)^{th}$ departure and the $k^{th}$ arrival on segment $j$.
The bounds $\underline{t}_j, \underline{s}_j, \underline{h}_j$ can be calculated based on the lower bounds above, applying relationships~(\ref{form1},\ref{form2},\ref{form3}).

Respectively, upper bounds used here, are denoted $\bar{w}_j, \bar{g}_j, \bar{h}_j$.
Let us denote by $\kappa$ the passenger load capacity per train. The upper bound $\bar{h}_j$ can be fixed as follows.
\begin{equation}
    \bar{h}_j := \kappa / (\sum_j \lambda^{in}_j), \quad \forall j,
    \label{eq-max-h}
\end{equation}
in such a way that all the passenger arrivals can theoretically be served.
Then, based on $\bar{h}_j$, the upper bounds $\bar{w}_j$ and $\bar{g}_j$ are fixed as follows.
\begin{align}
    & \bar{w}_j := x_j \bar{h}_j, \quad \forall j, \label{eq-max-w} \\
    & \bar{g}_j := \bar{h}_j - \underline{w}_j, \quad \forall j. \label{eq-max-g}
\end{align}

Let us finally note $\tilde{r}_j$ the nominal run time on segment~$j$.
$\tilde{r}_j$ is given by stability condition~(\ref{cond2}) below and depends on $\Delta{r}_j$.
Indeed $\Delta{r}_j$ influences the robustness and the frequency
of the metro system, i.e. it can either be fixed equal to $\Delta{w}_j$ to enhance the average frequency, or it can be chosen greater than $\Delta{w}_j$ to reinforce robustness towards perturbations.


\subsection{Review of the Max-plus linear traffic model~\cite{SFLG18a}}
\label{review1}

We shortly review here the Max-plus linear train dynamics for a metro loop line,
with demand-dependent dwell and run time controls~\cite{SFLG18a}. 
The dwell and run time control laws are fixed as follows.
\begin{align}
    & w^k_{j} = \min\{x_{j}h^k_{j} , \bar{w}_{j}\}, \label{dwell-acc-18} \\
    & r^k_{j} = \max\{\underline{r}_{j} , \tilde{r}_{j} - x_{j} (h^k_{j} - \underline{h}_{j})\}. \label{run-acc-18}
\end{align}

The authors of~\cite{SFLG18a} have shown that under the two following conditions
\begin{align}
    & h^1_{j} \leq \bar{h}^1_{j} = 1/(1-x_j) \bar{g}_{j}, \label{cond1} \\
    & \Delta r_{j} = \tilde{r}_j - \underline{r}_j \geq \Delta w_{j} = X_{j} \Delta g_{j}, \forall j, \label{cond2}
\end{align}
(\ref{dwell-acc-18}) and~(\ref{run-acc-18}) sum to
\begin{equation}\label{travel-acc-18}
  t_{j}^k = r_{j}^k + w_{j}^k = \tilde{r}_{j} + X_{j} \underline{g}_{j}.
\end{equation}

Then assuming that the $k^{th}$ departure from segment $j$ is realized as soon as the two
constraints on the travel and on the safe separation time are satisfied, the train dynamics are written as follows.
\begin{equation}\label{dynamics-acc-18}
   d^k_{j} = \max \left\{ \begin{array}{l}
                                d^{k-b_{j}}_{j-1} + \tilde{r}_{j} + X_{j} \underline{g}_{j}, \\
                                d^{k-\bar{b}_{j+1}}_{j+1} + \underline{s}_{j+1}.
                            \end{array} \right.
\end{equation}
It has been shown in~\cite{SFLG18a} that the train dynamics~(\ref{dynamics-acc-18}) admit a stationary regime
with a unique average growth rate, interpreted as the asymptotic average train time-headway $h$.
Moreover, the latter has been derived analytically by applying basic results of the Max-plus algebra theory~\cite{BCOQ92,CCGMQ98}, in function of $m$, the number of running trains on the metro line, $X$, the passenger demand,
and other parameters. We recall this principal result below.
\begin{equation}\label{av_h-acc-18}
    h(m, X) = \max \left\{\begin{array}{l}                           
			    \sum_j \frac{\tilde{r}_j + X_j \underline{g}_j}{m} ,\\~~\\
                            \max_j \left(\tilde{r}_j + X_j \underline{g}_j + \underline{s}_j\right) , \\~~\\
                            \sum_j \frac{\underline{s}_j}{n-m}
                           \end{array} \right.
\end{equation}

Under conditions~(\ref{cond1}) and~(\ref{cond2}), the dynamics are Max-plus linear and admits a stationary regime
with a unique average train time-headway $h(m,X)$ as shown in~(\ref{av_h-acc-18}).
These results are very interesting since we guarantee that under the control laws~(\ref{dwell-acc-18}) and~(\ref{run-acc-18})
on the train dwell and run times respectively, and under the two conditions~(\ref{cond1}) and~(\ref{cond2}), the
dwell times respond, i.e. are extended due to the passenger demand, whereas the run times are adapted in order to compensate
eventual extensions of the dwell times, to guarantee the overall train dynamics remain stable.

However, we know that the deviation in time of the train time-headway from its average value does not necessarily converge to zero,
but it can be periodic around zero.
Consequently, irregularities on the headway $h^k_{j}$ are not minimized over the $k^{th}$ departures.
This is due to the fact that the Max-plus linear system~(\ref{dynamics-acc-18}) may admit more than 
one Max-plus eigenvectors associated to its unique Max-plus eigenvalue.

\subsection{Review of the dynamic programming traffic model~\cite{FNHL16, FNHL17b, Far18}}
\label{review2}

We recall in this section some dynamic models proposed and studied in~\cite{FNHL16,FNHL17b,Far18}.
In those models, the train run times are assumed to be constant on all the segments of the metro line.
A first studied model fixed the train dwell times as follows.
\begin{equation}\label{unstable}
  w^k_j = \max\{\underline{w}_j , (\lambda_j / \alpha_j) g^k_j \}.
\end{equation}
Formula~(\ref{unstable}) models the fact that the train dwell times respect lower bounds $\underline{w}_j$ at every platform,
and are eventually extended to embark all the passengers accumulated during the close-in time $g^k_j$.
It has been recalled in~\cite{FNHL16} that the train dynamics with such dwell times is unstable.
The model~(\ref{unstable}) justifies the natural instability of the train dynamics, which has already
been shown in~\cite{Bres91}.

In~\cite{FNHL16, FNHL17b, Far18}, a new control law for the train dwell times at platforms has been proposed.
The dwell times have the following form.
\begin{equation}
\label{dwell-arxiv}
	w^k_j = \max\{\underline{w}_j , \bar{w}_j - \delta_j g^k_j \},
\end{equation}
where $\delta_j$ are parameters inversely proportional to $\lambda_j$, in such a way that $w^k_j$ will be
directly proportional to $\lambda_j$. An important fact here is that $w^k_j$ is inversely proportional 
to $g^k_j$, i.e. for a long $g_j$, the dwell time at the upcoming platforms is shortened.
This is advantageous to minimize headway irregularities without degrading the frequency, as it would be the
case for a holding control policy.
However, this is contra-intuitive from a passenger point of view, since, a retard leads to an accumulation
of passengers on the platform, i.e. trains entering the station with a long interval to the preceding train,
have to cope with a higher demand.

The authors of~\cite{FNHL16} have shown that with dwell time control~(\ref{dwell-arxiv}), and under the condition
$0\leq \delta_j\leq 1, \forall j$, the train dynamics
admits an asymptotic regime with a unique average growth rate interpreted here as the average train time-headway. 
Therefore the train dynamics are stable. Moreover, they are interpreted as the dynamic programming system of 
a stochastic optimal control problem of a Markov chain. 

\section{Traffic model with headway variance minimization}

On the one side the Max-plus linear traffic model of section~\ref{review1} from~\cite{SFLG18a} guarantees 
stability and controls both dwell and run times
to take into account the passenger travel demand. However, it does not harmonize train departure time intervals on the line.

On the other side, the dynamic programming model of section~\ref{review2} from~\cite{FNHL16, FNHL17b, Far18} guarantees stability but controls 
only the dwell times at platforms without controlling the run times. However, this model is better in term of harmonization
of train departure time intervals, comparing to the Max-plus one.

The new model we propose here extends the Max-plus model~(\ref{dynamics-acc-18}) of section~\ref{review1} to a dynamic programming
model as the model of section~\ref{review2}. Consequently, our new model benefits from the advantages of both models.
\begin{itemize}
\item Accounting for passenger demand with (partially) demand-dependent dwell times and run time control,
\item Minimizing train time-headway variances on the platforms via a dwell time control without degrading the average frequency.
\end{itemize}
%
%

We consider here the same train run time control law~(\ref{run-acc-18}), and propose a new train dwell time control law which modifies~(\ref{dwell-acc-18}) as follows.
%
\begin{align}
    w^k_j = \min \left\{(1-\gamma_j) x_j h^k_j , \bar{w}_j \right\},
    \label{dwell-eq}
\end{align}
with $0 \leq \gamma_j \leq 1$.

We notice that by activating the control ($0 < \gamma_j \leq 1$) in case of a perturbation on the train headway, we limit excessively long dwell times, which would have been a direct consequence of a long headway.

Let us use the following additional notations.
\begin{itemize}
  \item $\Delta w_j := \bar{w}_j - \underline{w}_j$.
  \item $\Delta g_j := \bar{g}_j - \underline{g}_j$.
  \item $\Delta r_j := \tilde{r}_j - \underline{r}_j$.
\end{itemize}
It is then easy to check that $\Delta w_j = 1/(1-x_j) \Delta g_j$.
\begin{proposition}~\label{proposition_1}
   If for all $j$, $h^1_j \leq \bar{h}_j$ and $\Delta r_j \geq \Delta w_j$, 
   then $t^k_j = \tilde{r}_j + X_j \underline{g}_j - \gamma_j X_j g^k_j, \forall j$.
\end{proposition}

\proof The proof is by induction. It is similar to the one of Theorem~1 in~\cite{SFLG18a}.
\endproof

The train dynamics is then written as follows.
\begin{equation}\label{dynamics}
    d^k_j = \max \left\{
	\begin{array}{l}
	    (1 - \delta_j) d^{k-b_j}_{j-1} + \delta_j d^{k-1}_j + (1-\delta_j) (\tilde{r} + X_j \underline{g}_j), \\~~\\
	    d^{k-\bar{b}_{j+1}}_{j+1} + \underline{s}_{j+1}.
	\end{array} \right.
\end{equation}
with $\delta_j := (\gamma_j x_j) / (1 + \gamma_j x_j)$.

Let us notice that the dynamics~(\ref{dynamics}) extends~(\ref{dynamics-acc-18}), since in the case where $\delta_j = 0$ (i.e. $\gamma_j = 0$),
(\ref{dynamics}) coincides with~(\ref{dynamics-acc-18}). This is a direct consequence of the fact that the train dwell time control law~(\ref{dwell-eq})
extends~(\ref{dwell-acc-18}).


It is easy to see that if $m=0$ (zero trains) or if $m=n$ (the metro line is full of trains), then the dynamic system~(\ref{dynamics}) is fully implicit. 
Indeed, if $m=0$, then $b_j = 0, \forall j$, i.e. the first term of the maximum operator in~(\ref{dynamics}) is implicit for every $j$.
Similarly, if $m=n$, then $\bar{b}_j = 0, \forall j$, i.e. the second term of the maximum operator in~(\ref{dynamics}) is implicit for every $j$.
We notice that in both cases, $m=0$ and $m=n$, no train movement is possible.
On the other side it is not difficult to check that if $0 < m <n$, then the dynamic system is implicit but triangular, i.e.
there exists an order on $j$ of updating the variables $d^k_j$ in such a way that the system will be explicit.
In fact, this order corresponds to the one of the train movements on the metro line.
In the following, we consider only the case $0<m<n$. Therefore, the dynamic system~(\ref{dynamics}) admits an equivalent triangular system.

The train dynamics~(\ref{dynamics}) can be written as follows.
\begin{align}\label{markov_chain}
    d^k_j = \max_{u \in U} \left\{ \left(M^u d^{k-1}\right)_j + \left(N^u d^k\right)_j + c^u_j \right\}, \forall j,k,
\end{align}
where $M^u, u\in\mathcal U$ and $N^u, u\in \mathcal U$ are two families of square matrices, and $c^u, u\in\mathcal U$ is a family of column vectors.
Moreover, since $0\leq \delta_j < 1, \forall j$ by definition, then $M^u_{ij} \geq 0, \forall u,i,j$ and $N^u_{ij} \geq 0, \forall u,i,j$.
Furthermore, we have $\sum_j \left(M^u_{ij} + N^u_{ij}\right) = 1, \forall u,i$.

The equivalent triangular system of system~(\ref{dynamics}) can be written as follows.
\begin{align}\label{markov_chain2}
    d^k_j = \max_{u \in U} \left\{ \left(\tilde{M}^u d^{k-1}\right)_j + \tilde{c}^u_j \right\}, \forall j,k,
\end{align}
where $\tilde{M}^u, u\in\mathcal U$ is a family of square matrices, and $\tilde{c}^u, u\in\mathcal U$ is a family of column vectors,
which can be derived from $M^u, N^u$ and $c^u, u\in\mathcal U$.
Moreover, we still guarantee $\tilde{M}^u_{ij} \geq 0, \forall u,i,j$ and $\sum_j \tilde{M}^u_{ij} = 1, \forall u,i$.
By consequent, the system~(\ref{markov_chain2}) can be seen as a dynamic programming system of an optimal control problem of a Markov chain,
whose transition matrices are $\tilde{M}^u, u\in\mathcal U$ associated to every control action $u\in\mathcal U$, and whose 
associated rewards are $\tilde{c}^u, u\in\mathcal U$.

\begin{theorem}~\label{thm-eigenvalue}
    If $0 < m < n$, then the dynamics~(\ref{dynamics}) admits a stationary regime, with a unique asymptotic average
    growth vector (independent of the initial state vector $d^0$) whose components are all equal to $h$, which is
    interpreted here as the asymptotic average train time-headway.    
\end{theorem}

\proof
We give only a sketch of the proof here. 
Since the dynamics~(\ref{dynamics}) are interpreted as the dynamic programming system of a stochastic optimal control problem
of a Markov chain,
it is sufficient to prove that the Markov chain in question, which is acyclic since $0 < m < n$,
is irreducible for every control strategy.
This is not obvious from the dynamics~(\ref{dynamics}). We need to show it on the equivalent triangular system~(\ref{markov_chain2}).
An alternative proof is the one of Theorem~5.1 in~\cite{FNHL16}.
\endproof

Theorem~\ref{thm-eigenvalue} does not give an analytic formula for the asymptotic average train time-headway $h$.
However, it guarantees its existence and its uniqueness. Therefore, $h$ can be approximated by numerical simulation
based on the value iteration, as follows.
\begin{align}
    h \approx d^K_j/K, \forall j, \text{ for a large } K.
\end{align}

Let us now go back to the objective of the extended model we propose here, i.e. our model~(\ref{dynamics}) 
extending the model~(\ref{dynamics-acc-18}) of~\cite{SFLG18a}.
We give here the main ideas explaining why the control we proposed permits to minimize the asymptotic average train time-headway variance.
Rigorous arguments with proved theorems will be given in our forthcoming articles.

As mentioned above, if $\gamma_j = 0, \forall j$, the dynamics~(\ref{dynamics}) are Max-plus linear, and are equivalent to~(\ref{dynamics-acc-18}).
In this case, if we write the dynamics~(\ref{dynamics}) or equivalently~(\ref{dynamics-acc-18}) under the triangular form~(\ref{markov_chain2}), then
the associated matrices $\tilde{M}^u$ are boolean circulant matrices. The latter may then have eigenvalue~$1$ with multiplicity bigger than~$1$.
Therefore, although the average growth rate of the dynamics is unique and independent of the initial state $d^0$, the asymptotic
state $d^k$ (up to an additive constant) may depend on the initial state $d^0$.

However, in case of the dynamics~(\ref{dynamics}) with $0 < \gamma_j \leq 1$ for some~$j$, the application of $\gamma_j$  will force the dynamics
to converge a stationary regime where the activated matrices $\tilde{M}^u$ are the ones
having eigenvalue~$1$ as a simple one, i.e. with multiplicity~$1$. In this case, $d^k$ will converge, up to an additive constant,
independent of the initial state $d^0$. Moreover, the asymptotic state $d$ will correspond to the case where the train time-headway
is the same at all the platforms. 

We will show in the next section by means of numerical simulation, that in this case,
the train dynamics converge to a state of traffic where the train departures are harmonized.
In other words, while the train time-headway converges to its asymptotic value $h$, the variance in time of the train time headways,
i.e. its deviation with respect to its average value $h$, converges to zero.  

With regard to dwell time equation~(\ref{dwell-eq}), we notice that in case $\gamma_j > 0$ where the control is applied, we do not only ensure train headway convergence at platforms, but are also sure not to degrade the train frequency. Indeed, with this, dwell times are slightly shortened to allow minimization of variances on the train time-headway. This means, in case of a big interval, train dwell times are no longer extended excessively, but are expanded just a little bit. On networks with a very high demand, as for the RATP network in Paris, this is highly advantageous to a holding control, because it allows a passenger load harmonization while maintaining a high train frequency.

\section{Simulation results}

\begin{table*}[thbp]
\centering
\vspace{10pt}
\caption{Simulation results for a metro loop line of RATP, Paris. In the left column, the train trajectories over all segments $j$ for $20$ trains,
in the right column the headways between two consecutive trains on the line at the end of the simulation, for $20, 21, 22$ trains.
In the first line without control ($\gamma_j = 0$), i.e. Max-plus linear traffic model
with demand-dependent dwell and run times, in the second line with static control
($\gamma_j = 0.1$), in the third line with dynamic control ($\gamma_j = 0.5 - (0.5/K)k$) for all $j$.}

\begin{tabular}{|c|c|}
  \hline
  $\gamma_j = 0$ & $\gamma_j = 0$ \\
  \hline
   & \\
  \includegraphics[scale=0.23]{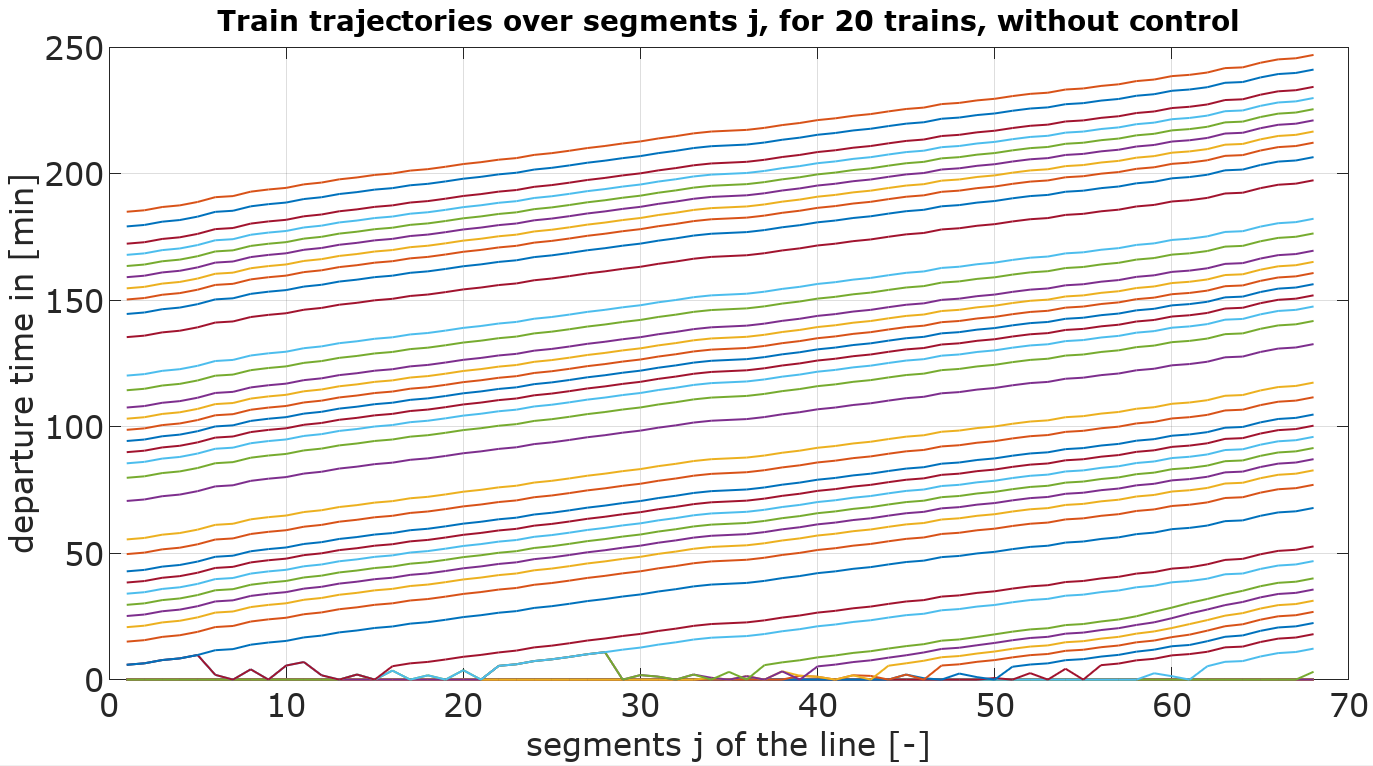} & \includegraphics[scale=0.23]{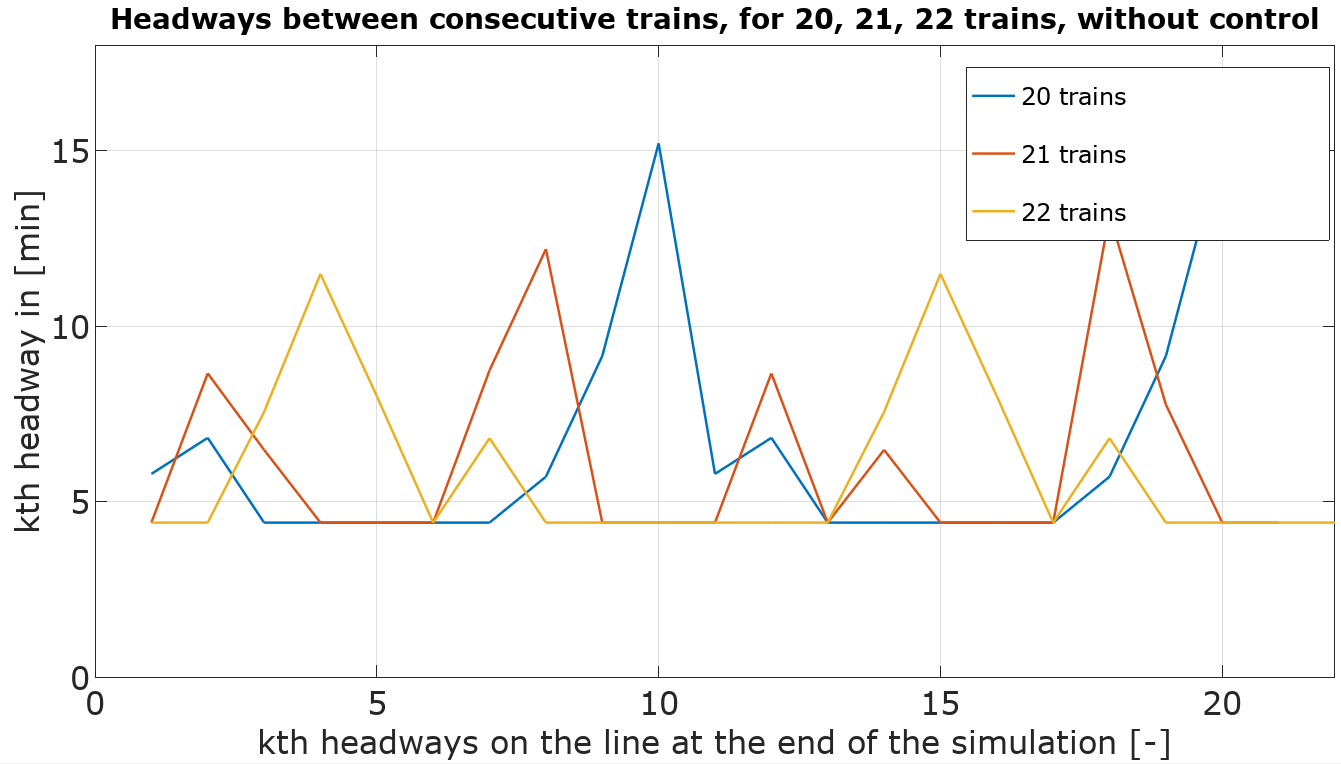} \\
  \hline 
  $\gamma_j = 0.1$ & $\gamma_j = 0.1$ \\
   \hline
   & \\
  \includegraphics[scale=0.23]{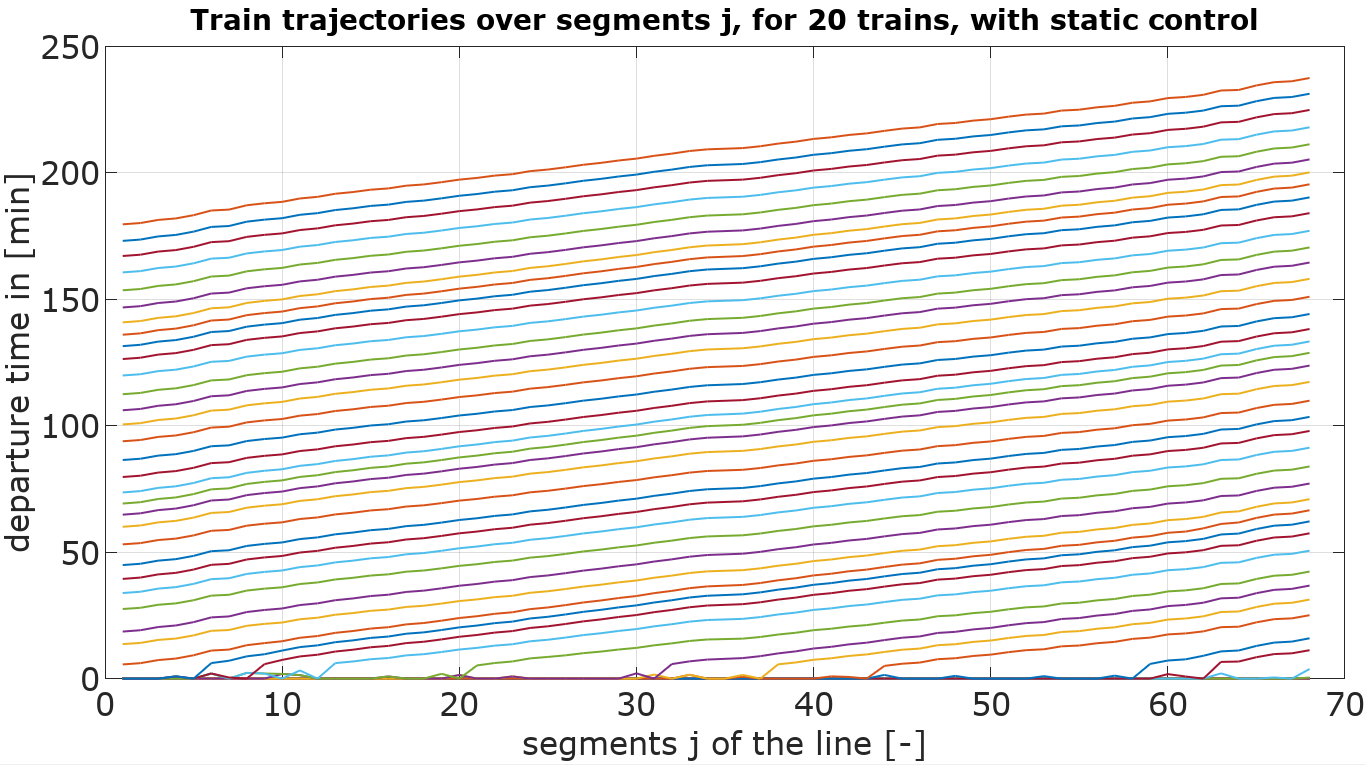} & \includegraphics[scale=0.23]{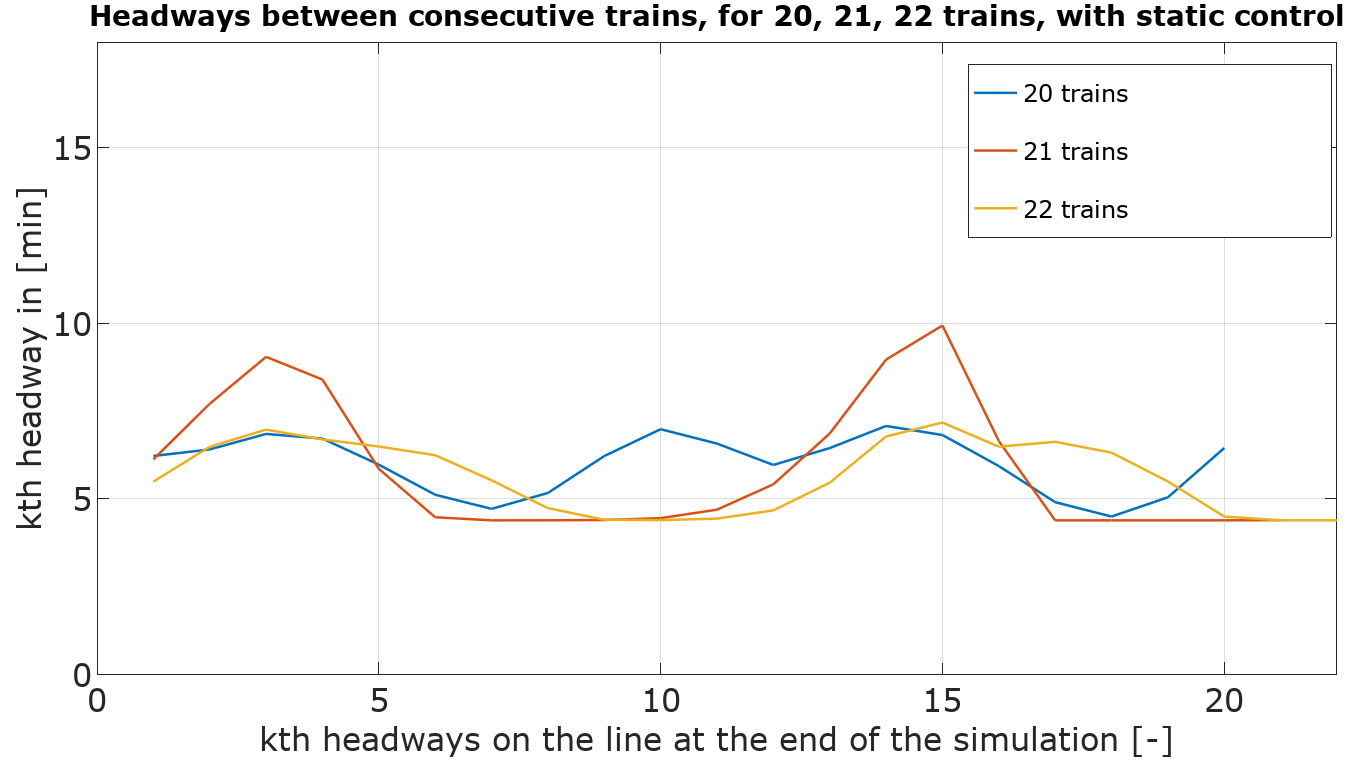} \\
  \hline 
  $\gamma_j = 0.5 -(0.5/K) k $ & $\gamma_j = 0.5 -(0.5/K) k$ \\
   \hline
   & \\
  \includegraphics[scale=0.23]{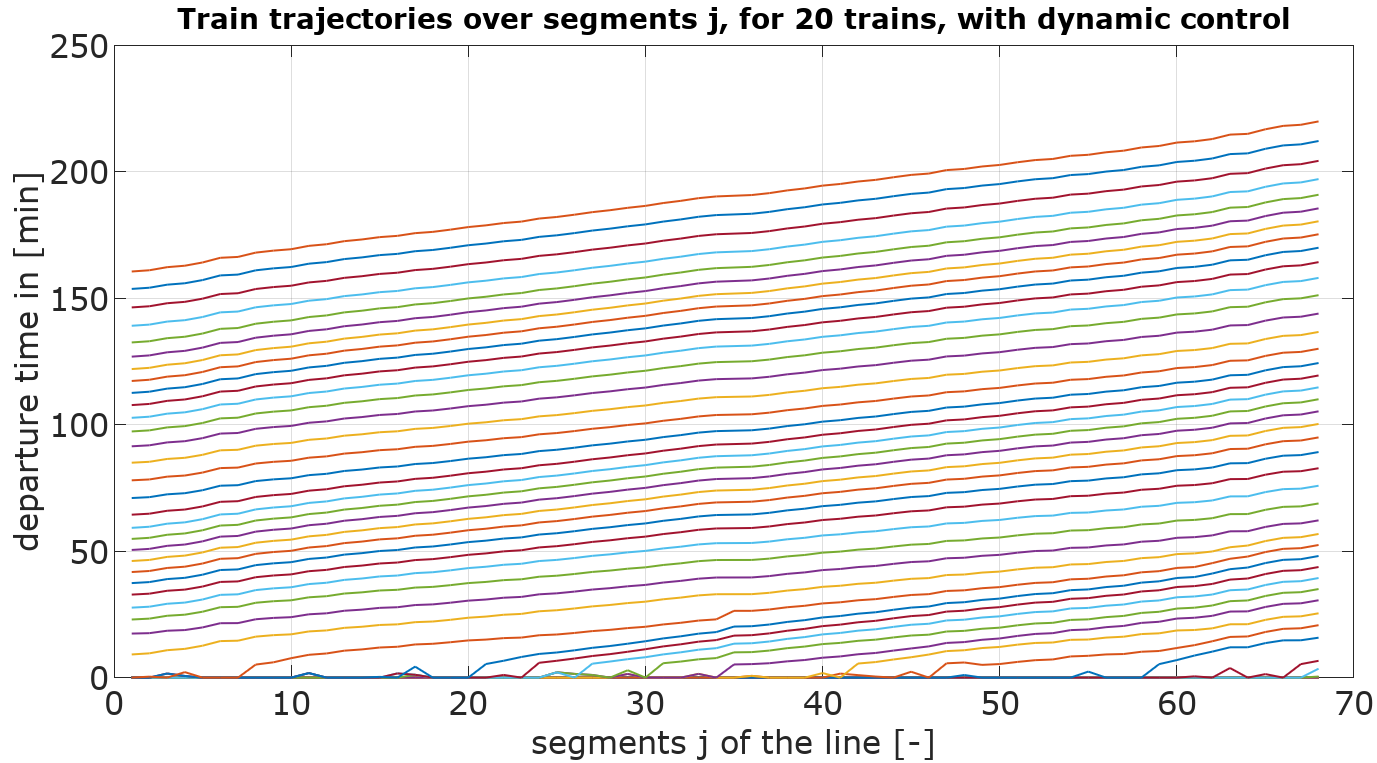} & \includegraphics[scale=0.23]{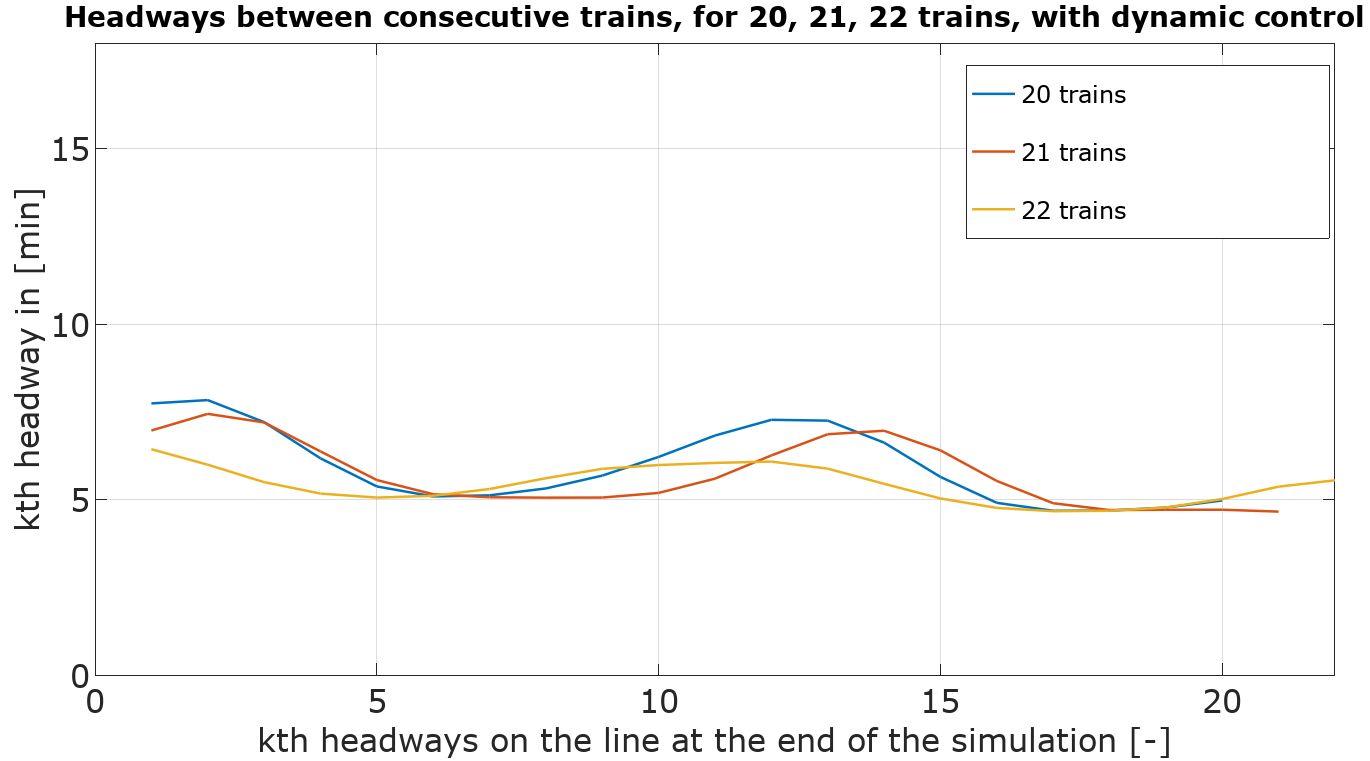} \\
  \hline
\end{tabular}
\label{tab}
\end{table*}


In Tab.~\ref{tab} we present some simulation results for different values of parameter $\gamma_j$.
First, we depict the results for the Max-plus dynamics, i.e. $\gamma_j = 0$.
Second, we fix $\gamma_j = 0.1$, which is close to the Max-plus dynamics.
Third, $\gamma_j$ degrades linearly from $0.5$ to $0$, over the simulation horizon $K$.
The simulation horizon $K$ is fixed accordingly to the time span, after which the train positions should be harmonized. 
It has been chosen to $K = 80$ departures here.
In the left column, we show the train trajectories, i.e. the departure times over all the segments on the metro loop line.
In the right column, we give the last observed headways $h(\gamma_j) = d(\gamma_j)^K_j - d(\gamma_j)^{K-1}_j$ over all the
segments $j$ on the line, indeed a good way to analyze the efficiency of the proposed control strategy, aiming to minimize headway variance.

The first row presents the results for $\gamma_j = 0$, i.e. without control where the dynamics~(\ref{dynamics}) are
a Max-plus linear system, as in~\cite{SFLG18a}. It can easily be seen, that, depending on the initial train positions,
we observe some areas with more trains on the line, and others with less, i.e. the headways are irregular. Moreover, this
irregularity persists over the simulation. In fact, it has been shown in~\cite{SFLG18a} that the train time-headways along the line
do not converge. They rather reach a periodic regime which makes their asymptotic average convergent. 
We notice that the train time-headways do not diverge although dwell times are modeled depending on the headway,
which remain constant with the applied run time control.
The graph giving the last observed headway shows clearly that important irregularities on the headway persist without control, 
depending on the number of trains, ranging from less than 5 minutes up to over 15 minutes.

The second row depicts the results for a static control where $\gamma_j = 0.1$ for all $j$ which include a platform.
In the left graph we see some irregularities at the beginning of the simulation, i.e. at the bottom, which tend to be minimized over time.
In fact, with regard to the final headways in the right graph, it can be seen that headways clearly converge along the line. However, due to the low $\gamma_j = 0.1$, larger initial variances cannot be minimized fast enough, i.e. they still can be observed at the end of the simulation. Variances finally range between less than 5 minutes and 10 minutes, i.e. a maximum variance of around 5 minutes.

In the last row, we give the results for a dynamic linear control where $\gamma_j = 0.5 -(0.5/K) k$ for all $j$ with a platform. The graph in the left column shows that some headway irregularities at the beginning are quickly being minimized and headways converge. In the right column, the graph giving the final headways observed along the line, shows that headways converge more quickly than for the static control. For the same simulation horizon as in the static case, the final headway variance does not exceed 2.5 minutes, i.e. headways range between slightly less than 5 minutes and around 7.5 minutes.
Finally, the dynamic control is advantageous because it applies stronger controls right after a perturbation in the transient regime, but returns to a $\gamma_j = 0$ once the perturbation has been absorbed, i.e. once the stationary regime has been reached. From a passenger point of view, this is preferable since for a $\gamma_j = 0$, dwell times respond once more fully to the time needed for passenger de- and embarking.

\section{Conclusion and future works}
Our series of traffic models for metro lines, ranging from a standard model with static run and dwell times, over an enhanced version with demand-dependent dwell times and a run time control, has been completed here by a real-time control with minimization of the train time-headway variance.
The simulation results have underpinned the effectiveness of the here proposed traffic control.
The figures clearly show headway convergence for the here proposed control. 
Moreover, they strongly suggest that a dynamic control, where a stronger control is applied right after a perturbation, followed by a fade-out time, leads to better results, i.e. faster headway variance minimization than a light but permanent control.

To complete this control model series, this paper will be followed by an application of this control policy to the traffic model for a metro line with a junction.
At this point, we assume to have a complete toolbox of traffic models which will allow an on-site implementation on the RATP network in Paris.
Possible applications are metro loop lines as well as metro lines with junctions, where the here presented model ensures headway regularity 
while maintaining a desired train frequency, with passenger demand-dependent dwell times and run time control.
There is always room for further improvement, e.g. for a passenger stock model in the trains and on the platforms, as well as an extension 
towards more complicated operational concepts, i.e. including intermediate terminus stations and non-stopping policy at selected stations.


\begin{thebibliography}{00}
\bibitem{BCOQ92} F. Baccelli, G. Cohen, G. J. Olsder, J.-P. Quadrat, "Synchronization and linearity: an algebra for discrete event systems," John Wiley and Sons, 1992.
\bibitem{Bres91} V. Van Breusegem, G. Campion, G. Bastin. "Traffic Modeling and State Feedback Control for Metro Lines," in IEEE Transactions on Automatic Control, Volume 36, Issue 7, Pages 770-784, 1991.
\bibitem{CCGMQ98} J. Cochet-Terrasson, G. Cohen, S. Gaubert, M. McGettrick, J.-P. Quadrat, "Numerical computation of spectral elements in maxplus algebra," IFAC Conference on System Structure and Control, Nantes, 1998.
\bibitem{FNHL16} N. Farhi, C. Nguyen Van Puh, H. Haj-Salem, J.-P. Lebacque, "Traffic modeling and real-time control for metro lines. arXiv:1604.04593, 2016.
\bibitem{FNHL17a} N. Farhi, C. Nguyen Van Puh, H. Haj-Salem, J.-P. Lebacque, "Traffic modeling and real-time control for metro lines. Part I - A Max-plus algebra model explaining the traffic phases of the train dynamics," IEEE American Control Conference, Seattle, 2017.
\bibitem{FNHL17b} N. Farhi, C. Nguyen Van Puh, H. Haj-Salem, J.-P. Lebacque, "Traffic modeling and real-time control for metro lines. Part II - The effect of passengers demand on the traffic phases," IEEE American Control Conference, Seattle, 2017.
\bibitem{Far18} N. Farhi, "Physical models and control of the train dynamics in a metro line without junction," in IEEE Transactions on Control Systems Technology, Issue 99, Pages 1-9, 2018. 
\bibitem{Fer05} A. Fernandez, A. P. Cucala, B. Vitoriano, F. de Cuadra, "Predicitive traffic regulation for metro loop lines based on quadratic programming," in Proceedings of the Institution of Mechanical Engineers, Part F: Journal of Rail and Rapid Transit, Volume  220, Pages 79-89, 2005.
\bibitem{Gov07} Rob M.P. Goverde, "Railway timetable stability analysis using 
max-plus system theory," in Transportation Research Part B, Volume 41, Issue 2, 179-201, 2007.
\bibitem{SFCLG17} F. Schanzenbacher, N. Farhi, Z. Christoforou, F. Leurent, G. Gabriel, "A discrete event traffic model explaining the traffic phases of the train dynamics in a metro line system with a junction," IEEE Conference on Decision and Control, Melbourne, 2017.
\bibitem{SFLG18a} F. Schanzenbacher, N.  Farhi, F. Leurent, G. Gabriel, "A discrete event traffic model explaining the traffic phases of the train dynamics on a linear metro line with demand-dependent control," IEEE American Control Conference, Milwaukee, 2018.
\bibitem{SCF16} F. Schanzenbacher, R. Chevrier, N. Farhi, "Train flow optimization: A model predictive and quadratic programming approach (Fluidification du trafic Transilien : approche pr\'{e}dictive et optimisation quadratique)," Conference ROADEF, Compi\`{e}gne, 2016.
\bibitem{LDYG17a} S. Li, M. M. Dessouky, L. Yang, Z. Gao, "Joint optimal train regulation and passenger flow control strategy for high-frequency metro lines," in Transportation Research Part B: Methodological, Volume 99, Pages 113-137, 2017.
\bibitem{LYG17b} S. Li, L. Yang, Z. Gao, "Optimal switched control design for automatic train regulation of metro lines with time-varying passengers arrival flow," in Transportation Research Part C, Volume 86, Pages 425-440, 2017.



\end{thebibliography}
\end{document}